\documentclass[manmat-mod,twoside]{svjour-mod}

\usepackage{amsfonts,amssymb,amscd,xy,graphicx}
\usepackage[leqno]{amsmath}
\usepackage{mathrsfs}
\def\bm#1{\mathpalette\bmstyle{#1}}
\def\bmstyle#1#2{\mbox{\boldmath$#1#2$}}

\newcommand{\homo}[2]{\mathrm {Hom}(#1,#2)}
\newcommand{\uhomo}[2]{\underline{\mathrm {Hom}}(#1,#2)}
\newcommand{\ext}[2]{{\mathrm {Ext }}^{1}(#1,#2)}
\newcommand{\extn}[2]{{\mathrm {Ext }}_{\mathbb Z/n\mathbb Z}^{1}(#1,#2)}
\newcommand{\biext}[3]{{\mathrm {Biext}}^1(#1,#2;#3)}

\newcommand{\extfl}[2]{{\mathrm {Ext }}^{1}_{\mathrm {R^\fl}}(#1,#2)}

\newcommand{\extlog}[3]{{\mathrm {Ext}}_{{#1}^\mathrm{log}}(#2,#3)}

\newcommand{\Cal}[1]{{\mathcal {#1}}}
\newcommand{\nobd}{\nobreakdash}

\newcommand{\id}{\mathrm{id}}

\newcommand{\fl}{\mathrm{fl}}
\newcommand{\cl}{\mathrm{cl}}

\newcommand{\loga}{\mathrm{log}}

\newcommand{\un}{\mathrm{un}}
\newcommand{\mgr}[1]{{\mathbb G}_{m, #1}}
\newcommand{\agr}[1]{{\mathbb G}_{a, #1}}

\newcommand{\zz}{{\mathbb Z}}
\newcommand{\qq}{{\mathbb Q}}
\newcommand{\FF}{{\mathbb F}}
\newcommand{\NN}{{\mathbb N}}
\newcommand{\MM}{{\mathbb M}}
\newcommand{\PP}{{\mathbb P}}

\newcommand{\spec}[1]{\mathrm {Spec}\left(#1 \right)}

\newcommand{\pre}[1]{{}_{#1}}
\renewcommand{\H}{\mathrm{H}} 

\newcommand{\add}{{\mathrm{add}}}
{\bf}{\it}
{\bf}{\it}
{\bf}{\it}
{\bf}{\rm}
{\bf}{\rm}
\newtheorem{thm-def}[lem]{Theorem(-Definition)}{\bf}{\rm}
{\bf}{\it}

\begin{document}
\input xy
\xyoption{all}
\title{Monodromy of logarithmic Barsotti-Tate groups attached to $1$\nobd-motives}
\author{
A. Bertapelle$^*$, M. Candilera, V. Cristante
}
 \institute{Dipartimento di Matematica Pura ed Applicata,
via Belzoni 7, I-35131 Padova\\
{e-mail: bertapel@math.unipd.it, candilera@math.unipd.it, cristanv@math.unipd.it}\\
$^*$ Partially supported by Progetto Giovani Ricercatori CPDG021784 Universit\`a di Padova
\\   }
\date{July 2, 2003}
\vskip 20pt

\maketitle

\vskip 30pt

\section*{Introduction}
 Let  $R$ be a complete discrete valuation ring with perfect
residue field $k$ of positive characteristic $p$ and field of
fractions $K$.
In this paper we consider a
$K$\nobd-1\nobd-motive $M_K$ as in \cite{[Ra]} and its associated
Barsotti-Tate group. This last does not in general extend to a
Barsotti-Tate group over $R$. However, with some assumptions, it
extends to a logarithmic Barsotti-Tate group over $R$. This
follows from \cite{[Ra]} and Kato's results on finite logarithmic
group schemes. Once chosen a uniformizing parameter $\pi$ of $R$,
any logarithmic Barsotti-Tate group over $R$ is described by two
data $(G,N)$ where $G$ is a classical Barsotti-Tate group over $R$
and $N$ is a homomorphism of classical Barsotti-Tate groups.
Moreover, if $R=W(k)$, $N$ induces a $W(k)$\nobd-homorphism
$\bm{\Cal N}\colon \MM(G_k)\to \MM(G_k)$ on Dieudonn\'e modules such
that $F\bm{\Cal N}V=\bm{\Cal N}$ and $\bm{\Cal N}^2=0$. In the
first part of the paper we recall these constructions and we show
how to relate $N$ with the ``geometric monodromy'' introduced by
Raynaud. In the second part of the paper we give an explicit
description of $\bm{\Cal N}$ in terms of additive
extensions and  integrals.
In the last part of the paper we
describe how to recover  the logarithmic Barsotti-Tate group
attached to a $1$\nobd-motive from a concrete scheme endowed with
a suitable logarithmic structure.

\section{$1$-motives}

\begin{definition}
\label{def.motive} Let $S$ be a scheme. An \emph{$S$\nobd-$1$\nobd-motive}
$M=[u\colon Y\to G]$ is a two term complex (in degree -1, 0)
 of commutative group schemes over $S$ such that:
\begin{itemize}\addtolength{\itemsep}{0.3\baselineskip}
\item $Y$ is an $S$\nobd-group scheme that locally for the \'etale
topology on $S$ is isomorphic to a constant group of type $\zz^r$,

\item $G$ is an $S$\nobd-group scheme extension of an abelian
scheme $A$ over $S$ by a torus $T$,

\item $u$ is an $S$\nobd-homomorphism $Y\to G$.
\end{itemize}

\end{definition}
Morphisms of $S$\nobd-$1$\nobd-motives are usual morphisms of complexes.
\pagebreak[3]

\begin{definition}\label{def.1mot}
Let $M_K$ be a $K$-$1$\nobd-motive. One says that $M_K$
\begin{enumerate}\addtolength{\itemsep}{0.5\baselineskip}

\item has \emph{good reduction} if $M_K$ extends to a
$1$\nobd-motive
    over $R$, i.e. if

\begin{itemize}\addtolength{\itemsep}{0.2\baselineskip}
  \item $Y_K$ is not ramified over $R$,
  \item $T_K$ has good reduction over $R$,
  \item $A_K$ has good reduction over $R$,
  \item $u_K$ extends to a homomorphism $u\colon Y\to G$.
\end{itemize}

(Hence $G_K$ extends to a semi-abelian $R$\nobd-group scheme $G$.)
 \item has \emph{semistable reduction} if

\begin{itemize}\addtolength{\itemsep}{0.2\baselineskip}
  \item $Y_K$ is not ramified over $R$,
  \item $T_K$ has good reduction over $R$,
  \item $A_K$ has semistable reduction over $R$.
\end{itemize}

(Hence $G_K$ extends to a smooth $R$\nobd-group scheme with semi\nobd-abelian
special fibre.\footnote{Cf. \cite{[Ra]} \S 4.})

\item has \emph{potentially semistable} (resp. \emph{good})
\emph{reduction} if it acquires semistable (resp. good) reduction
after a finite extension of $K$.

\item is \emph{strict} if $G_K$ has potentially good reduction.
\end{enumerate}

\end{definition}
Observe that any $K$\nobd-$1$\nobd-motive has potentially
semistable reduction. However, even if we allow base change, the
morphism $u_K$ does not in general extend over $R$. A simple
example is the Tate curve $u_K\colon \zz\to \mgr K$ with
$u_K(1)=\pi$ the uniformizing element.
 It has semistable reduction but no good reduction.

In the following we will consider only $K$\nobd-$1$\nobd-motives
or $R$\nobd-$1$\nobd-motives. For more details see Raynaud's paper
\cite{[Ra]}.

\subsection{The Barsotti-Tate group attached to a $K$-$1$-motive}

Let  $n$ be any positive integer and denote by $\pre n H$ the
kernel of $n$\nobd-multiplication on a group $H$.
For any $K$\nobd-$1$\nobd-motive
$M_K=[u_K\colon Y_K\to G_K]$ one can construct an exact sequence of
finite $n$\nobd-torsion group schemes over $K$:
\begin{equation}\label{eq.eta}
 \eta(n,u_K)\colon ~~~~~0\to \pre n G_K\to \pre n M_K\to Y_K/nY_K\to 0
\end{equation}
where $\pre n M_K$ is the cokernel of the homomorphism
\begin{equation*}
\displaystyle Y_K\overset{(-n,-u_K)}{\longrightarrow }
Y_K\times_{G_K}G_K ;
\end{equation*}
here the fibre product is taken with respect to $u_K$ on $Y_K$ and
the $n$\nobd-multiplication on $G_K$. As explained in \cite{[Ra]}
3.1, $\pre n M_K$ is the $\H^{-1}(C(M,n))$ with $ C(M,n)$ the cone
of the $n$\nobd-multiplication on the $1$\nobd-motive $M_K$, i.e.
\[\begin{array}{cccccl}C(M,n)\colon & Y_K&\longrightarrow
&Y_K\oplus G_K &\longrightarrow &G_K\\
                       & y&\mapsto& (-nx,-u_K(y))&&\\
                        & &&(y,g)&\mapsto &u_K(y)-ng
\end{array} \]
in degree $-2$, $-1$, $0$.

\begin{definition}
The $p$\nobd-divisible group or Barsotti-Tate group of the
$K$\nobd-$1$\nobd-motive $M_K$ is $\underset{\to}{\lim} (\pre {p^m} M_K)$.
\end{definition}
In the previous notations we have then an exact sequence of
BT-groups:
\[0\to \underset{\to}{\lim} (\pre {p^m} G_K)\to
\underset{\to}{\lim} (\pre {p^m} M_K)
\to \underset{\to}{\lim} ( Y_K/p^m Y_K)\to 0 .
\]
It is clear that if $M_K$ has good reduction then
$\underset{\to}{\lim} (\pre {p^r} M_K)$  extends to a
BT\nobd-group over $R$. We want to understand what
happens in the general case. We state now a result that we will
need later.


\begin{lemma}\label{lem.pullback}
Let notations be as above.
\begin{enumerate}
\item Consider the following diagram obtained via push\nobd-out by
$u_K$:
\[
\xymatrix{  0 \ar[r] & Y_K\ar[d]^{u_K} \ar[r]^{-n}& Y_K\ar[d] \ar[r]&
         Y_K/n Y_K\ar[d] \ar[r] & 0  \\
               0 \ar[r] & G_K \ar[r]& Y_K\amalg_{Y_K} G_K \ar[r]&
         Y_K/n Y_K \ar[r] & 0
}\] The short exact sequence  $ \eta(n,u_K)$ in (\ref{eq.eta}) is
isomorphic to the sequence of kernels for the
$n$\nobd-multiplication of the lower sequence.
\item Consider the
following diagram obtained via pull-back by
 $u_K$
\[\xymatrix{  0 \ar[r] & \pre n G_K \ar[d] \ar[r]& Y_K\times_{G_K} G_K
     \ar[d] \ar[r]&
         Y_K\ar[d]^{u_K} \ar[r] & 0  \\
               0 \ar[r] & \pre n G_K \ar[r]& G_K \ar[r]^n&
         G_K \ar[r] & 0
}\]
The short exact sequence  $ \eta(n,u_K)$ in (\ref{eq.eta}) is isomorphic to
 sequence of cokernels for the $n$\nobd-multiplication
 of the upper sequence.
\end{enumerate}
\end{lemma}

Raynaud shows in \cite{[Ra]} that to any $K$\nobd-$1$\nobd-motive
it is possible to associate in a
canonical way  a
$K$\nobd-$1$\nobd-motive with potentially good reduction $M_K'$
having the same BT\nobd-groups. His construction makes use of
rigid analytic methods. As a consequence, working with BT-groups
attached to a $K$\nobd-$1$\nobd-motive, \emph{one can always
assume the $K$\nobd-$1$\nobd-motive to be strict}. We will do so
in the sequel.


\subsection{Geometric monodromy}
Given a strict $K$\nobd-$1$\nobd-motive, the failure of good reduction is
controlled by a pairing, the so-called geometric monodromy.
To define it we need to recall some facts on the Poincar\'e bundle.

\begin{remark}\label{rem.1}
Let $M_K=[u_K\colon Y_K\to G_K]$ be a $K$\nobd-$1$\nobd-motive
with $Y_K^*$ be the group of characters of the torus part $T_K$ of
$G_K$ and $A_K$ the abelian  variety $G_K/T_K$. It is
known\footnote{See for
example \cite{[De]} 10.2.14 and \cite{[Ch]} II, 2.3.3.}
that to give a $1$\nobd-motive as above is equivalent to giving
morphisms $h_K\colon Y_K\to A_K$, $h_K^*\colon Y_K^*\to A_K^*$
(with $A_K^*$ the dual variety of $A_K$) and a trivialization
$s_K\colon Y_K\times Y_K^*\to {\cal P}_K$ of the pull-back
via $h_K\times h^*_K$ of the biextension ${\cal P}_K$.
Suppose that $G_K$ has good reduction. Then both $A_K$ and the dual abelian
variety $A_K^*$ have good reduction and the Poincar\'e bundle
${\cal P}_K$ extends to a biextension ${\cal P}$ in
$\biext{A}{A^*}{\mgr R}$ on N\'eron models. Also $h_K, h^*_K$
extend to morphisms $h,h^*$ over $R$ and the pull-back of ${\cal
P}$ via $h\times h^*$ provides a biextension ${\cal P}^Y $ in
$\biext{Y}{Y^*}{\mgr R}$. Moreover this is trivial on generic
fibres because of the existence of the trivialization $s_K$.
\end{remark}


\begin{definition}[\cite{[Ra]} \S4.3]\label{def.monodromy}
Let  $M_K=[u_K\colon Y_K\to G_K]$ be a strict  $K$\nobd-$1$\nobd-motive and
$Y_K^*$ the group of characters of $T_K$. The \emph{geometric
monodromy} of $M_K$ is a morphism
\begin{equation}\label{eq.mu}
\mu\colon Y_K \otimes Y_K^*\to \qq
\end{equation}
defined as follows:
\begin{enumerate}\addtolength{\itemsep}{0.5\baselineskip}
\item Suppose that $G_K$ has good reduction. Then ${\cal
P}^Y\in\biext{Y}{Y^*}{\mgr R}$ is trivial on generic fibres (see remark
\ref{rem.1}). Hence the biextension ${\cal P}^Y$ is the pull-back
of
\[ 0\to \mgr R\to {\cal G}\to i_*\zz\to 0
\]
 via a unique\footnote{Notations are those in \cite{[SGA7]} VIII; we used
that $\biext{Y}{Y^*}{{\cal G}}\cong\biext{Y_K}{Y^*_K}{\mgr K}$.}
$\mu_0\in \homo{Y_K\otimes Y^*_K}{\zz}= \homo{Y\otimes
Y^*}{i_*\zz}={\mathrm{ Biext}}^0({Y}, {Y^*}; {i_*\zz})$. One sets
$\mu=\mu_0$. \item In the general situation, $G_K$ reaches good
reduction after a Galois extension $K'$ of $K$.  Now the monodromy
on $K'$ is compatible with Galois action and can be descended to a
$\mu$ as in (\ref{eq.mu}).
\end{enumerate}
\end{definition}
Observe that  $\qq$ has to be thought of as the group of values of the
valuation of the algebraic closure of $K$ with $\zz$  the group
of values assumed on $K$.

   Let $K^\un$ be the maximal unramified extension of $K$, $v\colon
(K^\un)^*\to \zz$ the valuation and $R^\un$ its valuation ring.
Observe that in the hypothesis of $i)$ there is a valuation
$v_{{\cal P}}$ on ${\cal P}_K(K^\un)$ and that $\mu_0=v_{{\cal
P}}\circ s_K$ holds. Moreover if the abelian part is trivial, then
${\cal P}_K=\mgr K$, $h_K$ and $h_K^*$ are the structure morphisms
and $s_K\colon Y_K\otimes Y_K^*\to \mgr K$ is the usual pairing
 $(y,y^*)\mapsto y^*(y)$. Hence $\mu_0(y,y^*)=v(y^*(y))$.
 These results can be
generalized. See also 4.6/6 in \cite{[Ra]}.


\begin{lemma}\label{lem.monodromy}
Let notations be as above. Suppose that $G_K$ has good reduction
over a finite field extension $L$ of $K$. Then the
geometric monodromy pairing $\mu$ coincides with the pairing
\begin{equation*}
 Y_K\otimes Y_K^*\to \qq, \;\;\;(y,y^*)\mapsto \frac{1}{e}v_L(y^*(t))
\end{equation*}
where $e$ is the index of ramification of $L/K$, $v_L$ is the
extension of the valuation of $L$ to $L^\un$,
$t\in T_K(L^\un)$ is any point whose image
in the component group $\zz^d$ of $T_{ L^\un}$ (as well as of $G_{ L^\un}$)
coincides with the image of $u_K(y)$ and $y^*(t)\in ( L^\un)^*$ for any
$y^*\colon T_{ L^\un}\to \mgr { L^\un}$.
\end{lemma}
\proof  We may reduce to the case $L=K$. As explained above
$\mu_0(y,y^*)$ is obtained via the valuation on ${\cal P}_K(K^\un)$.
Consider the push-out
\begin{eqnarray}
\xymatrix{
0 \ar[r] & T_K \ar[d]^{y^*} \ar[r] & G_K \ar[d]^{g_{y*}}
             \ar[r] & A_K  \ar@2{-}[d]   \ar[r] & 0 \\
0 \ar[r] & \mgr K  \ar[r] & G_{y^*}
             \ar[r] & A_K    \ar[r] & 0
}
\end{eqnarray}
where  $G_{y^*}\cong{\cal P}_{K|A_K\times h^*(y^*)}$. The
valuation $v_{\cal P}$ on ${\cal P}_K(K^\un)$ restricts to a $v_{y^*}$ on
$G_{y^*}(K^\un)$ and $\mu_0(y,y^*)=v_{y^*}(g_{y*}\circ u_K(y))$ holds.
Let now $t\in T_K(K^\un)$ be a point having the same image as
$u_K(y)$ in the component group $\zz^d$. Then one has is
$v_{y^*}(g_{y*}\circ u_K(y))=v(y^*(t)$ for any character $y^*$ and
\mbox{$u_K(y)-t\in G(R^\un)$}.
To conclude it is now sufficient to observe that
$v_{y^*}\circ g_{y*} $
is zero on $G(R^\un)$ and that \mbox{$v(y^*(t))=v_{y^*}(g_{y*}t)$}.
 \qed


\subsection{Devissage}
Once having realized that the defect of good reduction is controlled by
the geometric monodromy, Raynaud explains, under the hypothesis that
the geometric monodromy takes integer values, how to decompose
a strict $1$\nobd-motive into the sum of two $1$\nobd-motives, the first
having potentially good reduction and the second codifying the
monodromy.


\begin{theorem}\label{thm.devissage}
Let $M_K=[u_K\colon Y_K\to G_K ]$ be a strict $K$\nobd-$1$\nobd-motive such
that the geometric monodromy $\mu$ factors through $\zz$.
Then for any choice of a uniformizing parameter $\pi$ of $R$
there is a canonical decomposition
\begin{equation}\label{eq.devissage}
u_K=u_{K,\pi}^1+u_{K,\pi}^2
\end{equation}
where $u_{K,\pi}^2$ factors through the torus part $T_K$ and is
given by the formula
\begin{eqnarray}\label{eq.u2}
u_{K,\pi}^2 \colon  Y_K &  \to & T_K= \uhomo{Y_K^*}{\mgr K}  \overset{\iota}{\to} G_K\nonumber  \\
 y &\mapsto & (y^*\mapsto \pi^{ \mu(y,y^*) })
\end{eqnarray}
while $u_{K,\pi}^1$ has potentially good reduction.
\end{theorem}
\proof  \cite{[Ra]}, 4.5.1 \qed

\begin{remark}\label{rem.devissage}
If both $G_K$ and $Y_K$ have good reduction, the geometric
monodromy factors through $\mathbb Z$ and the $1$\nobd-motive
$u_{K,\pi}^1$ in the previous decomposition
has good reduction.
\end{remark}

\begin{remark} Let notations be as above and
$\rho\colon G_K\to A_K$ be the morphism of $G_K$ in its abelian
quotient. It is clear that $\rho\circ u_K=\rho\circ u_{K,\pi}^1$
has potentially good reduction. Moreover the push-out of
$\eta(n,u_K)$ with respect to $\rho_n\colon \pre n G_K\to \pre n
A_K$ (the restriction of $\rho$ to the kernels of $n$\nobd-multiplication)
is $\eta(n,\rho\circ u_K)$. More
precisely we have:
\begin{equation}\label{dia.HVM}
\xymatrix{ &   & \pre n T_K \ar[d]^w \ar@2{-}[r]& \pre n T_K \ar[d]^\tau \\
 \eta(n,u_K)\colon &0 \ar[r] &\pre n G_K \ar[d]^{\rho_n} \ar[r] &
              \pre n M_K \ar[d]^g \ar[r]^h & Y_K/nY_K \ar@2{-}[d] \ar[r] & 0\\
\eta(n,\rho\circ u_K)\colon &0 \ar[r] &\pre n A_K  \ar[r] &
              \pre n M_K^A \ar[r]^f & Y_K/nY_K \ar[r] & 0
}
\end{equation}
\end{remark}
\bigskip

Suppose that the hypothesis of Theorem~\ref{thm.devissage} holds.
 We wish to compare the  sequence of  finite $n$\nobd-torsion
 group schemes $\eta(n,u_K)$ in (\ref{eq.eta}) associated to $u_K$ with the
sequences associated to  $u_{K,\pi}^1$ and $u_{K,\pi}^2$.


\begin{lemma}\label{lem.baersum}
Let $M_K\colon [u_K\colon Y_K\to G_K ]$ be a strict
$K$\nobd-$1$\nobd-motive such that the geometric monodromy factors
through $\zz$. Let $ u_K=u_{K,\pi}^1+u_{K,\pi}^2$ be the
decomposition of Theorem~\ref{thm.devissage}. Then $\eta(n,u_K)$
is isomorphic to $\eta(n,u_{K,\pi}^1)+\eta(n,u_{K,\pi}^2)$
 where $+$ denotes Baer's sum.
\end{lemma}
\proof Consider the homomorphism
\begin{equation*}
\homo{Y_K}{G_K}\overset{\partial}{\longrightarrow }\ext{Y_K}{\pre
n G_K}\cong\extn{Y_K/nY_K}{\pre n G_K}
\end{equation*}
that associates to a $1$\nobd-motive $u_K$ the pull\nobd-back of
$0\to \pre n G_K\to G_K\stackrel{n}{\to} G_K\to 0$ by $u_K$, resp.
the sequence of cokernels of such pull-back. Here the subscript
$\zz/n\zz$ stands  for extensions in the category of
$\zz/n\zz$\nobd-modules. We have already seen in
Lemma~\ref{lem.pullback}~ii) that the isomorphism class of
$\eta(n,u_K)$ is $\partial(u_K)$.  The result follows from the
fact that $\partial$ is a homomorphism. \qed

\begin{lemma}\label{lem.vanishingtwo}
Given a $K$-$1$\nobd-motive $u_K\colon \zz^r\to \mgr K^d$, (the
isomorphism class of) the sequence $\eta(n,u_K)$ extends over $R$
if and only if $u_{K,\pi}^2$ is divisible by $n$, i.e. if and only
if $\eta(n,u_{K,\pi}^2)$ is isomorphic to the trivial sequence.
\end{lemma}
\proof
Recall the following diagram
\[
\xymatrix{\homo{\zz^r}{\mgr R^d} \ar[d] \ar[r]^n & \homo{\zz^r}{\mgr R^d}\ar[d] \ar[r]^\partial &
\ext{\zz^r}{ {\bm \mu}_{n,R}^d } \ar[d] \ar[r] & 0 \\
\homo{\zz^r}{\mgr K^d} \ar[r]^n \ar[d] & \homo{\zz^r}{\mgr K^d}  \ar[r]^\partial\ar[d]  &
\ext{\zz^r}{{\bm \mu}_{n,K}^d} \ar[r] & 0 \\
\zz^{rd} \ar[r]^n & \zz^{rd} &
}
\]
where as above $\partial$ is obtained by pull-back. Suppose now that
 $\eta(n,u_K)=\partial(u_K)$ extends over $R$. The same is true for
 $\eta(n,u_{K,\pi}^1)$ because $u_{K,\pi}^1$ has good reduction and
 hence also $\eta(n,u_{K,\pi}^2)$
extends over $R$. Let $w_K$ be a $1$-motive with good reduction
such that $\eta(n,w_K)=\eta(n,u_{K,\pi}^2)$. Define
 $w_K^\prime := w_K-u_{K,\pi}^2$; observe that this is also the Raynaud decomposition
of $w_K^\prime$ and that $\partial(w_K')=0$. Hence there exists a
$1$\nobd-motive $w_K^{\prime\prime}$ such that $n\cdot
w_K^{\prime\prime}=w_K^\prime$. Let
\mbox{$w_K^{\prime\prime}=w_{K,\pi}^1+w_{K,\pi}^2 $} be Raynaud's
decomposition. It is clear that the $n$-multiplication preserves
Raynaud's decompositions and hence $n\cdot w_{K,\pi}^1=w_K$ while
$n\cdot w_{K,\pi}^2=-u_{K,\pi}^2$. Hence $u_{K,\pi}^2$ is
divisible by $n$. \qed


\subsubsection{Geometric monodromy \'a la Kato.}
In the following we will work with group schemes as sheaves
 of $\mathbb Z$\nobd-modules on the flat site. We wish to understand better
 the $K$-$1$\nobd-motive $u^2_{K,\pi}$ in (\ref{eq.devissage})  in order
to compare Kato's monodromy and Raynaud's monodromy.

Given an \'etale group scheme $N_K$ isomorphic to $\zz^r$ over an
algebraic closure of $K$, denote by  $N^\vee_K$ the \'etale group
scheme $\uhomo{N_K}{\zz}$ and by $N_K^{\vee D}$ its  Cartier dual.
We have $\zz^{\vee D}=\mgr K$ and $N_K^{\vee D}=N_K\otimes_\zz \mgr
K$. The  geometric monodromy $\mu\colon Y_K\otimes_\zz Y^*_K\to
\zz$ of $u_K$, provides  a morphism
\begin{eqnarray}\label{eq.nu}
\nu\colon Y_K &\longrightarrow  &\uhomo{Y_K^*}{\mathbb Z} =
\colon (Y_K^*)^\vee \\
y&\mapsto & \mu(y,-) \nonumber
\end{eqnarray}
and hence a morphism of tori
\begin{equation}\label{eq.nutori}
\nu\otimes \id
\colon \quad Y_K\otimes_\zz \mgr K
\to( Y_K^*)^\vee\otimes_\zz \mgr K=(Y_K^*)^D =T_K.
\end{equation}
 Let $H_K(1)$ denote the Cartier dual of the
Pontrjagin dual\footnote{The Pontrjagin dual of $H_K$ is
$\uhomo{H_K}{\qq/\zz}$.} for any finite \'etale $K$\nobd-group
scheme $H_K$. Then ${\bm \mu}_n=\zz/n\zz(1)$ and if $n$ kills
$H_K$ one has $H_K(1)=H_K\otimes_{\zz/n\zz} {\bm \mu}_n$. Hence we
can introduce a ``monodromy" homomorphism of level $n$
\begin{equation}\label{eq.nun}
\nu_n\colon  Y_K/nY_K(1)= Y_K\otimes_\zz {\bm \mu}_n
\longrightarrow  ( Y_K^*)^\vee\otimes_\zz {\bm \mu}_n =\pre n T_K
~~\left(\hookrightarrow \pre n G_K\right)
\end{equation}
as the restriction of   $\nu\otimes \id$ to the $n$\nobd-torsion
subgroups. It was defined in \cite{[Ra]} 4.6.

 Consider now the
 Tate $1$\nobd-motive
\[\bm \pi\colon \zz\to \mgr K, ~ 1\mapsto \pi.
\]
It is clear from (\ref{eq.u2}) that  $u_{K,\pi}^2$
has the following factorization
\begin{equation}\label{eq.vmotive}
\begin{array}{ccccccc}
 Y_K&~~\overset{{\bm \pi}_Y:=\id_{Y_K}\otimes{\bm \pi}}{\longrightarrow }~~
&Y_K\otimes_\zz \mgr K &~~ \overset{\nu\otimes
\id}{\longrightarrow }~~
& (Y_K^*)^\vee\otimes_\zz \mgr K=T_K& \stackrel{\iota}{\longrightarrow }& G_K\\
y&\mapsto& y\otimes \pi&\mapsto & \mu(y,-)\otimes \pi=\pi^{\mu(y,-) } & &
\end{array}
\end{equation}
where  $\nu$ was  defined in
(\ref{eq.nu}).

\begin{lemma}
Let $G_K$ have good reduction $G$, $\iota\colon T_K\to G_K$ be the
torus part and suppose that $Y_K$ extends to an \'etale group $Y$
over $R$. Then the homomorphism
\begin{eqnarray}\label{eq.hom}
 \homo{Y}{ G}\times
  \homo{Y_K\otimes \mgr K}{ T_K} &{\longrightarrow }&
\homo{Y_K}{G_K}\nonumber \\
(u^1,w_K)\hspace{3cm}& \mapsto& ~~u^1_K +\iota\circ w_K\circ
{\bm \pi}_{Y}
\end{eqnarray}
is an isomorphism.
\end{lemma}
The homomorphism ${\bm \pi}_Y$ was defined in (\ref{eq.vmotive}).
\proof The surjectivity is clear because any homomorphism $u_K$ on
the right hand side represents a $K$\nobd-$1$\nobd-motive and we have just
seen how to decompose $u_K$ as sum of a $1$\nobd-motive
 $u_{K,\pi}^1$ and a $1$\nobd-motive
$u_{K,\pi}^2=\iota\circ (\nu\otimes \id)\circ {\bm \pi}_{Y_K}$
with $\nu$ the monodromy homomorphism. In our hypothesis
$u_{K,\pi}^1$ extends to a $R$\nobd-$1$\nobd-motive $u_\pi^1$.
Hence $u_K$ is the image of the pair $(u_\pi^1,\nu\otimes \id)$.
For the injectivity: the first group on the left injects into
$\homo{Y_K}{G_K}$, so we are reduced to showing that given a
$w_K\in\homo{Y_K\otimes \mgr K}{ T_K}$ the
$K$\nobd-$1$\nobd-motive $\iota\circ w_K\circ {\bm \pi}_Y$ has
good reduction if and only if $w_K$ is trivial. Denote by
 $\mu(w_K)$ the  pairing corresponding to $w_K$ via the canonical isomorphisms
$$
\homo{Y_K\otimes \mgr K}{
T_K}=\homo{Y_K^*}{ Y_K^{\vee}}=\homo{Y_K\otimes Y_K^*}{\zz}.
$$
The image of $(0,w_K)$ via the map in (\ref{eq.hom})
is a $K$\nobd-$1$\nobd-motive; let  $\mu(w_K)$
denote its geometric monodromy. Such $K$\nobd-$1$\nobd-motive has
good reduction if and only if $\mu(w_K)=0$ and this last occurs if and
only if $w_K=0$. \qed
\medskip

We restrict again to the consideration of the $1$\nobd-motive
 $\bm \pi\colon \zz\to \mgr K, ~ 1\mapsto \pi$. Observe that it satisfies the
hypothesis of remark \ref{rem.devissage} and that in this case
$u^1_{K,\pi}$ is trivial. Let denote by
\begin{equation}\label{eq.thetak} \theta^\pi_{n,K}\colon 0\to {\bm \mu}_n\to \pre n E_K \to
\mathbb Z/n\mathbb Z\to 0
\end{equation}
the short exact sequence  $\eta(n,{\bm \pi})$. The following
results will be used in Theorem \ref{thm.comparison} to compare
 Raynaud's monodromy and Kato's monodromy.


\begin{theorem}\label{thm.push}
Let  $n$ be a positive integer and $\theta^\pi_{n,K}$  the short
exact sequence just defined. Suppose that $u_K$ is a strict
$K$\nobd-$1$\nobd-motive as in Theorem~\ref{thm.devissage} with
$u_K=u_{K,\pi}^1+u_{K,\pi}^2$ its Raynaud decomposition.
\begin{enumerate}\addtolength{\itemsep}{0.5\baselineskip}
\item The geometric monodromy of $u_K$  and the one of
$u^2_{K,\pi}$ coincide. \item The short exact sequence
$\eta(n,u_{K,\pi}^2)$ associated to $M^2_K=[u_{K,\pi}^2\colon
Y_K\to T_K]$ is isomorphic to the push-out  via $\nu_n$ of the
sequence
\[\theta^\pi_{n,K}\otimes_{\zz/n\zz} Y_K/nY_K\colon~ 0\to
\bm \mu_n\otimes_{\zz/n\zz} Y_K/nY_K\to \pre n
E_K\otimes_{\zz/n\zz} Y_K/nY_K\to Y_K/nY_K \to 0\]
\end{enumerate}
\end{theorem}
\proof The first fact follows immediately from the definition of
$u_{K,\pi}^2$. For the second assertion, consider the
factorization of $u_{K,\pi}^2$ described in  (\ref{eq.vmotive}).
It says that there is a commutative diagram
\[\xymatrix{
\eta(n, {\bm \pi}_Y)\colon & 0 \ar[r] & \pre n(Y^\vee_K)^D \ar[d]^{\nu_n} \ar[r]
     & \pre n M_K^\pi \ar[d] \ar[r]  & Y_K/nY_K \ar@2{-}[d] \ar[r] &0 \\
\eta(n,u_{K,\pi}^2)\colon  & 0 \ar[r] & \pre n T_K \ar[r] & \pre n M_K^2
 \ar[r] & Y_K/nY_K \ar[r] & 0
}\]
where $M_K^\pi=[{\bm \pi}_Y\colon Y_K\to Y_K\otimes \mgr K]$.
It is clear that
$\eta(n,{\bm \pi}_Y)=\theta^\pi_{n,K}\otimes Y_K/nY_K$ by definition of ${\bm \pi}_Y$.
Hence $\eta(n,u_{K,\pi}^2)$ is isomorphic to ${\nu_n}_*(\theta^\pi_{n,K}\otimes Y_K/nY_K)$, \emph{i.e} to
the push-out via $\nu_n$ of the
sequence $\theta^\pi_{n,K}\otimes_{\zz/n\zz} Y_K/nY_K$.
\qed


\begin{corollary}\label{cor.split}
Suppose furthermore that $Y_K$ and $G_K$ have
good reduction and consider the following homomorphism
\begin{equation}\label{eq.split}
\begin{array}{crclcc}
 \Psi\colon & \extn{\frac{Y}{nY}}{\pre n G}& \times
  &\homo{\frac{Y_K}{nY_K}(1)}{\pre n G_K}&{\longrightarrow }&
\extn{\frac{Y_K}{nY_K}}{\pre n G_K} \\
& & & & &
\\
&(\eta^1&,&h) &\mapsto& ~~\eta^1_K +
h_*(\theta^\pi_{n,K}\otimes_{\zz/n\zz} \frac{Y_K}{nY_K})
\end{array}
\end{equation}
where  $\eta^1_K$ means the restriction of $\eta^1$ to generic
fibres.
\begin{enumerate}\addtolength{\itemsep}{0.5\baselineskip}
\item If $u_K\colon Y_K\to G_K$ is a $K$\nobd-$1$\nobd-motive with
Raynaud decomposition $u_K=u_{K,\pi}^1+u_{K,\pi}^2$, then the
class of $\eta(n,u_{K})$ lies in the image of $\Psi$. More
precisely it corresponds to the pair $( \eta(n,u^1_\pi),\nu_n )$
where  $\nu_n$ is the ``monodromy'' homomorphism of level $n$ as
in (\ref{eq.nun}) and  $u^1_\pi$ is the $R$\nobd-$1$\nobd-motive
that extends $u_{K,\pi}^1$. \item If $Y_K\cong \zz^r$ and $G_K\cong \mgr
K^d$, then $\Psi$ is  an isomorphism.
\end{enumerate}
\end{corollary}
\proof  The first assertion is an immediate consequence of the
previous Theorem, part i). We restrict then to the case $\pre n
G_K= {\bm \mu}_n^d$ and $Y_K/nY_K=\zz^r/n\zz^r$.
 For the surjectivity it is sufficient to remark that any extension
class on the right is represented by a
 $\eta(n,u_{K})$ for a strict $K$\nobd-$1$\nobd-motive $u_K$
because of the vanishing of $\H^1(K,\mgr K^d)=\ext{\zz}{\mgr
K^d}$. For the injectivity: the group of extensions on the left
injects in the group of extensions on the right. It remains to
check that $h_*(\theta^\pi_{n,K}\otimes \zz^r/n\zz^r)$ extends
 over $R$ if and only if $h=0$. Now, $h\colon {\bm \mu}_n^r \to{\bm \mu}_n^d$
 extends to many homomorphisms $\tilde h_K\colon \mgr K^r\to \mgr K^d$.
Choose one of them and let $u_K\colon \zz^r \to\mgr K^d$ be
 $\tilde h_K\circ \bm\pi_{\zz^r}$
with $\bm\pi_{\zz^r}$ as in (\ref{eq.vmotive}). It is clear that
$h$ coincides with the monodromy homomorphism of level $n$ of the
$K$\nobd-$1$\nobd-motive $u_K$. By hypothesis $\eta(n,u_K)\cong
h_*(\theta^\pi_{n,K}\otimes~\zz^r/n\zz^r)$ extends over $R$. Hence
$u_K=u_{K,\pi}^1+u_{K,\pi}^2$ with $u_{K,\pi}^2$ divisible by $n$
(see Lemma \ref{lem.vanishingtwo}).
 This implies that the monodromy of $u_K$ that equals the monodromy of
$u_{K,\pi}^2$  is a multiple of $n$ and hence its monodromy homomorphism
of level $n$ is trivial, i.e. $h=0$.
 \qed
 \medskip

Theorem \ref{thm.push} ii) and its corollary are the only original
results of the first part of this paper. They become interesting
once one realizes that Kato proves an analogous result for extensions
of finite logarithmic group schemes (cf.~Theorem~\ref{thm.kato}).
The comparison of these two results makes it possible to extend $\pre
n M_K$ to a finite logarithmic group scheme over $R$.
\smallskip

We close this section by giving an example that should clarify all
the previous constructions.


\subsubsection{Tate's curve.} Let $n$ be a positive integer and
\begin{equation*}
u_K\colon \zz\to \mgr K; \quad 1\mapsto q=\epsilon \pi^{nr+s},
\quad 0\leq s\leq n-1,\quad 0\leq r, \quad \epsilon \in R^*
\end{equation*}
 an elliptic curve with split multiplicative reduction. The
canonical decomposition of Theorem~\ref{thm.devissage} provides
$u_{K,\pi}^2\colon \zz\to \mgr K,~ 1\mapsto \pi^{nr+s}$ and
\mbox{$u_{K,\pi}^1 \colon \zz\to \mgr K,~ 1 \mapsto\epsilon$}. The
geometric monodromy $\mu\colon \zz\otimes\zz\to \zz$ depends only
on $u_{K,\pi}^2$ and it sends $1\otimes 1$ to \mbox{$rn+s$}. The
``monodromy" homomorphism of level $n$, $\nu_n\colon
\zz/n\zz\otimes {\bm \mu}_n={\bm \mu}_n\to {\bm \mu}_n$, is the
$s$\nobd-multiplication. It is also clear that
$\eta(n,u_{K,\pi}^2)$ is isomorphic to $s\cdot\theta_{n,K}^\pi$.

\section{Finite logarithmic  group objects}

For the theory of logarithmic spaces we refer to \cite{[Il]} and \cite{[K3]}.
 We need also some definitions and results in \cite {[K1]}, \cite
{[K2]}.

Let $\pi$ be a fixed uniformizing element of $R$ and $\underline
T$ the spectrum of $R$ with the standard log structure given by
the chart $\NN \to R$, $1\mapsto \pi$. Denote by $T^\loga_\fl$ the
logarithmic flat site over $\underline T$. A finite
(representable) logarithmic group $G$ over $R$ is a sheaf of
abelian groups over $T^\loga_\fl$ that is represented by a fine
saturated log-scheme over $R$, log flat and of Kummer type over
$R$ so that its underlying scheme is finite over $R$. For an
example, consider Tate's elliptic curve  $E_K$ defined via ${\bm
\pi}\colon \zz\to \mgr K,~ 1\mapsto \pi$. Kato shows how to
extend $E_K$ to a group object $\underline E^\pi$ in the category
of valuative logarithmic spaces over $R$. This is explained by
Illusie  in \cite{[Il]} 3.1. The kernel of $n$\nobd-multiplication
on  $\underline E^\pi$, denoted by $\pre n(\underline E^\pi)$, is
obtained via log blow-ups from a  logarithmic space having
\begin{equation}\label{eq.finite}
\pre n E=\spec{ \oplus_{i=0}^{n-1}
\frac{R[x_i]}{(x_i^n-\pi^i)} }
\end{equation}
as underlying scheme. Moreover there is a short exact
sequence of finite logarithmic groups given by
\begin{equation}\label{eq.theta}
\theta_n^\pi\colon \quad 0\to \zz/n\zz(1)\to \pre n(\underline
E^\pi) \to \zz/n\zz\to 0.
\end{equation}
(cf. \cite{[Il]} 3.2.1.4) whose restriction to generic fibres is
the short exact sequence
$\theta^\pi_{n,K}$ that we used in
 Theorem~\ref{thm.push} (cf. \cite{[Il]} 3.2.1.4.).

 Let now  $F$ (resp. $H$) be
$n$\nobd-torsion finite (resp. finite \'etale) group schemes over $R$
endowed with the inverse image log structure. A  result of Kato
(cf. \cite{[K1]} p.~84) says that extensions
(of sheaves in $T^\loga_\fl$)
\begin{equation*}
\eta^\loga\colon \quad 0\to F\to G^\loga \to H\to 0
\end{equation*}
correspond bijectively (up to isomorphisms) to  pairs $(G^\cl, N)$ where
\begin{equation*}
\eta^\cl\colon \quad 0\to F\to G^\cl \to H\to 0
\end{equation*}
is a classical extension of group schemes over $R$ and
$N\colon H(1)\to F$ is a morphism of $R$\nobd-group
schemes where $ H(1)=H \otimes_{\zz/n\zz}{\bm \mu}_n$.
Moreover $\eta^\loga$ is the
Baer sum of $\eta^\cl$ and the push-out by $N$ of the extension
$\theta^\pi_n\otimes_{\zz/n\zz} H $.

\begin{theorem}[Kato]\label{thm.kato}
Let notations be as above. There is an isomorphism
\begin{eqnarray*}
  \extfl  {H}{F} \times
  \homo{H(1)}{F}&\overset{\sim}{\longrightarrow }&
\extlog{T} {H}{F} \\
(\eta^\cl,N)\hspace{1cm}&\mapsto& ~~\eta^\cl +N_*(\theta^\pi_n\otimes H)
\end{eqnarray*}
\end{theorem}
\proof cf. \cite{[K1]}.\qed
\smallskip

Observe that the statement of this theorem is similar to that of
Corollary~\ref{cor.split}.
We will explain in Theorem \ref{thm.comparison} that they are deeply
related. Before proceeding we need the following result.

\begin{lemma}\label{prop.comparison}
Let $M_K=[u_K\colon Y_K\to G_K ]$ be a $K$\nobd-$1$\nobd-motive and
 suppose that $G_K$ extends to a semiabelian $R$\nobd-scheme $G$
while $Y_K$ extends to an \'etale $R$\nobd-group scheme $Y$.
 Having fixed a uniformizing parameter $\pi$ of $R$,
let $u_K=u_{K,\pi}^1+u_{K,\pi}^2$ be Raynaud's decomposition
of Theorem \ref{thm.devissage}. Then:
\begin{enumerate}\addtolength{\itemsep}{0.5\baselineskip}
\item $u_{K,\pi}^1$
extends to a $R$\nobd-$1$\nobd-motive $u_{\pi}^1$;
 \item The monodromy homomorphism of level $n$ of $u_K$, i.e.
 $\nu_{n }\colon Y_K/nY_K(1)\to \pre n G_K$, extends to a homomorphism
$\nu_{n,R}\colon Y/nY (1)\to \pre n G $.
\end{enumerate}
\end{lemma}
\proof Both assertions are evident because
it follows from the hypothesis that the torus part
$T_K$ of $G_K$ extends to a torus over $R$; hence also its
group of characters $Y_K^*$ extends to an \'etale group over $R$,
say $Y^*$. This implies that the geometric monodromy takes values in $\zz$
and it extends to a biadditive map $Y\otimes Y^*\to \zz$ over
$R$. This last provides the homomorphism  $\nu_{n,R}$ we are looking for.
\qed \medskip

 We can now state the relation between Raynaud's geometric monodromy
of a $K$\nobd-$1$\nobd-mo\-tive and Kato's monodromy of its
logarithmic BT-group.

\begin{theorem}\label{thm.comparison} Let $M_K=[u_K\colon Y_K\to G_K ]$ be a
$K$\nobd-$1$\nobd-motive and suppose that $G_K$ extends to a
semiabelian $R$\nobd-scheme $G$ while $Y_K$ extends to an \'etale
$R$\nobd-group scheme $Y$. Let $\nu_{n,R}$ be the homomorphism
in Lemma \ref{prop.comparison}
with  $u^1_\pi\colon Y\to G$ the
$R$-$1$-motive that extends the $u^1_{K,\pi}$ of Raynaud's
decomposition.

The sequence $\eta(n,u_K)$ in
(\ref{eq.eta}) extends (up to isomorphisms) to a sequence of
finite logarithmic group schemes and precisely, in the notations
of Theorem \ref{thm.kato}, to the one associated to
$(\eta(n,u^1_\pi), \nu_{n,R} )$.

In particular, $\nu_n$ (i.e. Raynaud's monodromy
homomorphism of level $n$  associated to $u_K$) is Kato's
monodromy homomorphism $N$ restricted to generic fibres.
\end{theorem}
\proof
We know
from Lemma \ref{lem.baersum} and Theorem \ref{thm.push} that
$\eta(n,u_K)$ is isomorphic to
\[\eta(n,u_{K,\pi}^1)+(\nu_n)_*(\theta^\pi_{n,K}\otimes Y_K/nY_K).\]
Moreover, as $u_{K,\pi}^1$ extends to an $R$-$1$-motive $U_\pi^1$
(cf. Lemma \ref{prop.comparison})
also  $\eta(n,u_{K,\pi}^1)$ extends to a
sequence of classical group schemes $\eta(n,u_{\pi}^1)$; on the other hand
the sequence $\theta^\pi_{n,K}$ in (\ref{eq.thetak}) extends to
the sequence of logarithmic groups
that we denoted by $\theta^\pi_{n}$ in (\ref{eq.theta}). Hence
the sequence
of logarithmic groups
$\eta(n,u_{\pi}^1)+(\nu_{n,R})_*(\theta^\pi_n\otimes Y/nY)$ restricted
to generic fibres is (up to isomorphisms)
$\eta(n,u_K)$.
\qed

As all constructions above behave well  with respect to
inclusion homomorphisms \mbox{$\pre {p^m} M_K\to \pre {p^{m+1}} M_K$}
we can conclude that:

\begin{corollary}\label{cor.ext}
With hypothesis  as above, let  $M^1$ denote the $R$\nobd-$1$\nobd-motive
$[u_\pi^1\colon Y\to G]$. The BT-group of $u_K$,
$\underset{\to}{\lim} (\pre {p^m} M_K)$, extends to a logarithmic BT-group
$\underset{\to}{\lim}(\pre {p^m} M^1, \nu_{p^m,R})$ where
$\underset{\to}{\lim}\pre {p^m} M^1$ is the BT-group of $M^1_\pi$.
\end{corollary}

\begin{remark}
Let $F=(\zz/n\zz)^r$ and $H={\bm \mu}_n^d$. The decomposition in
Theorem \ref{thm.kato}  restricted to generic fibres coincides
with the isomorphism of Corollary~\ref{cor.split}. In particular
this is true working with the $K$-$1$-motive $E_K=[u_K\colon
\zz\to \mgr K]$  given by $1\mapsto q=\epsilon\pi^{nr+s}$,
\mbox{$0\leq s\leq n-1, \epsilon \in R^*$}. The kernel of
$n$\nobd-multiplication $\pre n E_K$ extends (up to isomorphism)
to a finite logarithmic group $(\pre n E^\cl,N)$ where $N$ is the
$s$\nobd-multiplication on ${\bm \mu}_n$ and $\pre n E^\cl$ is the
finite group scheme that lies in the middle of $\eta(n,u^1_{\pi})$
for $u_\pi^1\colon\zz\to \mgr R,~ 1\mapsto \epsilon$.
\end{remark}

\section{Monodromy on Dieudonn\'e modules}

Throughout this section $R=W(k)$ and
 $M_K=[u_K\colon Y_K\to G_K]$ will be a fixed
$K$-$1$-motive with $Y_K\cong \zz^r$.
We will suppose that $G_K$ has good reduction and its extension over $R$, $G$,
has split torus part of rank $d>0$.
We have seen in Corollary \ref{cor.ext} that under these hypotheses
the BT-group ${\bm M}_K(p):=\underset{\to}{\lim} (\pre {p^m} M_K)$
associated to $u_K$  extends to a
logarithmic BT-group, say ${{\bm M}(p)}^\loga=({\bm M}(p), N)$ where
${\bm M}(p)=\underset{\to}{\lim} (\pre {p^m} M^1_\pi)$
is a classical BT-group over $R$  and
\[N\colon Y\otimes_\zz \bm{\mu}_{p^\infty} \to
Y^{*\vee}\otimes \bm{\mu}_{p^\infty} \to\underset{\to}{\lim} (\pre {p^m} G)
\]
is the monodromy homomorphism:
the first morphism is $\nu_R\otimes \id$ where $\nu_R$
is the extension of $\nu\colon Y_K\to  Y^{*\vee}_K$ in
(\ref{eq.nu}) over $R$, which exists since
 $Y_K, Y_K^*$ are unramified.

\subsection{The identification of Dieudonn\'e modules of $\qq_p/\zz_p$ and
 $\bm{\mu}_{p^\infty}$}

For the theory of Dieudonn\'e modules we refer to \cite{[Fo]}:
If $G$ is a BT-group over $k$, its  Dieudonn\'e module is
$\MM(G)=\homo{G}{\widehat{CW}_k}$, where
Hom means homorphisms of $k$\nobd-formal groups. In particular,
$\MM((\qq_p/\zz_p)_k)=\zeta W(k)$ is a free $W(k)$\nobd-module of rank $1$
whose canonical generator $\zeta$ is the natural embedding of
$(\qq_p/\zz_p)_k$ in $\widehat{CW}_k$ once  $\qq_p/\zz_p$ is  identified with
$CW(\FF_p$); more precisely, $\zeta$ corresponds to the covector
 $y=(\dots,y_{-2},y_{-1})\in \widehat{CW}_k(k^{\qq_p/\zz_p})$
defined by $y_{-i}=\sum_{a \in \qq_p/\zz_p}a_{-i}f_a$, where
$f_a(b)=\delta_{ab}$ is the  Kronecker delta.

Let us recall that elements of  $\MM(G)$ can also be described as
 isomorphism classes of rigidified extensions of $G$ by $\agr k$
(as fppf sheaves,  cf.  \cite{[Kz]} \S 5.2 or \cite{[MM]} \S15).
Given a $\varphi \in \MM(G)$, the corresponding
 additive  extension  is obtained as the pull-back via
$\varphi$ of the extension
\[0\to \agr k\to \widehat{CW}_k \stackrel{V}{\to}  \widehat{CW}_k \to 0 \]
where $V$ is the Verschiebung of $\widehat{CW}_k$;
it will be denoted by $\varphi^\add_{p^\infty}$ . In particular,
one can prove that the extension
$$\zeta^\add_{p^\infty}:0\to \agr k\to F\to (\qq_p/\zz_p)_k\to0,$$
is isomorphic  to the push-out of
\begin{equation}\label{eq.zeta}
\zeta_{p^\infty}\colon 0\to
\zz\to \zz[1/p]\stackrel{f}{\to} \qq_p/\zz_p \to 0
\end{equation} via the canonical homorphism
$\zz \to \agr k$.
Let $\sigma\colon\qq_p/\zz_p\to \zz[1/p] $ be the section of
 $f$ such that
$0\leq\sigma(a)<1 $, for $a\in \qq_p/\zz_p$.
The factor set  of $\zeta_{p^\infty}$ (and hence of $\zeta^\add_{p^\infty}$)
corresponding to  the section $\sigma $ is then
\begin{equation}\label{eq.gamma}
\gamma\colon  (\qq_p/\zz_p) \times  (\qq_p/\zz_p) \to \zz ~(~ \to \agr R ),
\quad \left(a, b\right)\mapsto [\sigma(a)+\sigma(b)]\end{equation}
where square brackets mean integral  part.
The extension $\zeta^\add_{p^\infty}$ has  a canonical lifting
$(\zeta^\add_{p^\infty})_R$ to $R$
obtained as the push-out of $\zeta_{p^\infty}$ in (\ref{eq.zeta})
 (now as a sequence over $R$)
via the morphism $\zz\to \agr R$. The restriction of
$(\zeta^\add_{p^\infty})_R$ on generic fibres splits, and the map
\[h\colon \underset{\to}{\lim}\zz/p^i\zz=\qq_p/\zz_p \to  \agr K,
\quad a\mapsto \sigma (a)\]
 is the trivialisation. Let $K[X]$
be the affine algebra of $\agr K$ and let
$h^*_i: K[X]\to K^{\zz/p^i\zz}$, $i>0$,  be the
corresponding $K$-homomorphisms; one gains an element
\begin{equation}\label{eq.h}
h^*(X)=
\sum_{a\in \qq_p/\zz_p}\sigma(a)f_a\in K^{\qq_p/\zz_p}=\underset{\gets}{\lim}K^{\zz/p^i\zz}.
\end{equation}
If $\PP$ denotes the coproduct of $K^{\qq_p/\zz_p}$, then
\[\PP h^*(X)-1\widehat \otimes h^*(X)-h^*(X)\widehat \otimes1=
\sum_{a,b \in \qq_p/\zz_p}[\sigma (a)+\sigma (b)]f_a\widehat \otimes f_b\in
R^{\qq_p/\zz_p}\widehat \otimes R^{\qq_p/\zz_p};\]
this tells us that $h^*(X)$ is an integral of  second kind\footnote{We
recall that given the algebra $\cal A$ of a formal group
over $R$, an \emph{integral of
the second kind} of $\cal A$ is an element $f\in {\cal A}\hat\otimes_R K$
such that $df\in \Omega_{R}({\cal A})$ and
${\mathbb P}f-f\hat \otimes 1-1\hat\otimes f\in {\cal A}
\hat\otimes {\cal A} $, where $\mathbb P$ denotes the coproduct in $\cal A$;
we will denote by $I_2({\cal A})$ the $R$-module of integrals of the
second kind. Moreover if ${\mathbb P}f-f\hat \otimes 1-1\hat\otimes f=0$,
$f$ is called an \emph{integral of the first kind}.} of $ R^{\qq_p/\zz_p}$.

Moreover, the covector $y$, and hence $\zeta$, can be recovered from $h^*(X)$. In fact by  a direct  computation one can check that
$h^*(X)=\sum_{i=0}^ \infty p^{-i}\widehat y_{-i}^{p^i}$, where $\widehat y_{-i}\in R^{\qq_p/\zz_p}$ is a lifting of $ y_{-i}\in k^{\qq_p/\zz_p}$.
 \smallskip

Also the Dieudonn\'e module of $\bm \mu_{p^\infty}$ is a free $W(k)$\nobd-module
of rank 1. Let $R[[Y]]$ be the affine algebra of $\bm{\mu}_{p^\infty,R}$,
where $Y$ is the canonical parameter
and let $l(Y)\in \zz_{(p)}[\![Y]\!]$ be the  Artin--Hasse logarithm of $1+Y$,
\emph{ i.e.}
$$\text{exp}(-l(Y)-p^{-1}l(Y)^p-p^{-2}l(Y)^{p^2}-\dots )=1+Y, $$
then,
$\MM((\bm{\mu}_{p^\infty})_k)=\delta W(k)$, where $\delta=(\dots, l_0(Y_0),l_0(Y_0))$, and $l_0(Y_0)$ is the image of $l(Y)$ in the affine algebra of
$\bm \mu_{p^\infty,k}$.

Let us remark that
$-\text{log}(1+Y)= l(Y)+p^{-1}l(Y)^p+p^{-2}l(Y)^{p^2}+\dots$
is the integral of first kind of
$(\bm\mu_{p^\infty})_R$ obtained  by lifting $\delta$.

Let $(\qq_p/\zz_p)_k^\vee=\varinjlim (\zz/p^n\zz)_k^\vee$ be the Pontrjagin
dual of $(\qq_p/\zz_p)_k$; then $\bm{\mu}_{p^\infty,k}$ is the Cartier dual of
 $(\qq_p/\zz_p)_k^\vee$;
 as a consequence there exist two perfect pairings of $W(k)-$modules:
$$\langle -,-\rangle_C:\MM(\bm{\mu}_{p^\infty,k})\times  \MM((\qq_p/\zz_p)_k^\vee)\to W(k)\quad
\langle -,-\rangle_P:\MM((\qq_p/\zz_p)_k)\times  \MM((\qq_p/\zz_p)_k^\vee)\to W(k).$$
We will denote by
\begin{equation}\label{eq.id1}
\text{id}(1): \MM(\bm{\mu}_{p^\infty,k}) \to
\MM((\qq_p/\zz_p)_k)\end{equation} the $W(k)-$isomorphism such that
$$\langle m,n\rangle_C=\langle \text{id}(1)(m),n\rangle_P$$
for every $m\in \MM(\bm{\mu}_{p^\infty,k}) $ and
$n\in  \MM((\qq_p/\zz_p)_k^\vee)$. One can check that
$\text{id}(1)(\delta)=\zeta$, so that
\[F\circ \text{id}(1)\circ V=\text{id}(1).
\]


\subsection{Monodromy and additive extensions}

We want to define a $W(k)$\nobd-endomorphism of the Dieudonn\'e module
 $ \MM({\bm M}(p)_k)$  depending on the monodromy $N$.
(See also  \cite{[K1]} 5.2.2.)
We proceed as follows:
\[
\xymatrix{ \MM(Y^{*\vee}\otimes {\bm \mu}_{p^\infty,k})
\ar[r]^{\MM(N_k)} &
\MM(Y\otimes_\zz {\bm \mu}_{p^\infty,k}) \ar[r]^{\id_Y(1)}&
\MM(Y\otimes_\zz \qq_p/\zz_p)_k\ar[d]\\
\MM({\bm M}(p)_k) \ar[u]\ar@{..>}[rr]^{\bm{\Cal N}}& &
\MM({\bm M}(p)_k)  }
\]
where the vertical map on the left comes from the obvious
inclusion \[Y^{*\vee}\otimes \bm \mu_{p^n,k}= \underset{\to}{\lim}
(\pre {p^m} T_k) \to \underset{\to}{\lim} (\pre {p^m} M_{\pi,k}^1)\]
and the vertical map on the right is obtained from the projection
\[\underset{\to}{\lim} (\pre {p^m} M_{\pi,k}^1)\to Y\otimes_\zz
(\qq_p/\zz_p)_k.\]
 Recall now that both $Y$ and $Y^*$ are constant groups
and hence
\begin{equation}\label{eq.decomp1}
\MM(Y^{*\vee}\otimes \bm \mu_{p^\infty,k})=
Y^*\otimes \MM(\bm \mu_{p^\infty,k}), \quad
\MM(Y\otimes {\bm \mu}_{p^\infty,k})=Y^\vee \otimes_\zz
\MM({\bm \mu}_{p^\infty,k}).
\end{equation}
The morphism $\MM(N_k)$ is then
$$ \MM(N_k)=\nu^\vee\otimes \id ,$$
where $\nu^\vee\colon Y^{*} \to Y^\vee$ comes from the geometric
monodromy $\mu$ and is the transpose of $\nu$ in (\ref{eq.nu})
(on special fibres).
Using decompositions as in (\ref{eq.decomp1}) we can define
\begin{equation}\label{eq.id1Y}
\id_Y(1) := \id\otimes \id(1)\colon
Y^\vee\otimes \MM(\bm{\mu}_{p^\infty,k})\to Y^\vee \otimes \MM((\qq_p/\zz_p)_k)
 \end{equation}
where $\id(1)$ is the canonical identification explained in
(\ref{eq.id1}).
The map $\id_Y(1)$ is an isomorphism of
$W(k)$\nobd-modules such that
$F\circ \id_Y(1)\circ V=\id_Y(1)$ and
 the dotted arrow $\bm{\Cal N}$ turns out to be a
$W(k)$\nobd-homomorphism such
that ${\bm{\Cal N}}^2=0$ and $F{\bm{\Cal N}}V={\bm{\Cal N}}$.
\smallskip

We want now to describe how the composition $\id_Y(1)\circ \MM(N_k)$
works: Given an element $\chi\otimes \delta \in Y^*\otimes
\MM({\bm \mu}_{p^\infty,k})$, with $\delta$ the canonical generator
of $\MM({\bm \mu}_{p^\infty,k})$,
\[\left(\id_Y(1)\circ \MM(N_k)\right)(\chi\otimes \delta) =
\id_Y(1)\left( \nu^\vee(\chi)\otimes \delta\right)
=
\nu^\vee(\chi)\otimes \zeta .
\]
\begin{remark} The construction above could  also be done restricting to
kernels of $p^n$-multiplication $\pre {p^n} M$ and hence working
with the monodromy homomorphism of level $p^n$,
\[\nu_{p^n}\colon Y\otimes \zz/p^n\zz\to Y^*\otimes {\bm \mu}_{p^n}.\]
This is what Kato does in  \cite{[K1]} 5.2.2. Hence the construction
above is simply a way to summarize Kato's construction for all $p^n$.
\end{remark}

The element $\nu^\vee(\chi)\otimes \zeta\in  Y^\vee \otimes
\MM((\qq_p/\zz_p)_k)=
\MM(Y \otimes(\qq_p/\zz_p)_k)$ can be described in different ways; we give
here some examples without proofs:

\begin{enumerate}
\item If we interpret it as extension of $Y \otimes(\qq_p/\zz_p)_k$ by
the additive group $\agr k$, $\nu^\vee(\chi)\otimes \zeta$ is represented by
 the pull-back of the sequence $\zeta_{p^\infty}^\add$ with respect to
\[(\nu^\vee(\chi), \id) \colon
Y\otimes_\zz  \qq_p/\zz_p\to  \qq_p/\zz_p. \]
\item Given an extension of $\qq_p/\zz_p$ by $\agr k$, the push-out with
respect to the $m$\nobd-multiplication and the pull-back with
respect to the $m$\nobd-multiplication provide isomorphic
sequences. In a similar way one proves that the sequence
$(\nu^\vee(\chi), \id)^*\zeta_{p^\infty}^\add$
 is also
isomorphic to the push-out with respect to
$(\nu^\vee(\chi),\id)$ of the sequence
$Y \otimes\zeta_{p^\infty}^\add $.
\item
 Starting with the $1$\nobd-motive $\bm \pi\colon \zz\to \mgr K$,
consider the sequence
\[0\to \mgr K\to F_K\to \qq_p/\zz_p\to 0
\]
obtained first applying push-out $\bm \mu_{p^n}\to \mgr K$ to the sequence
$\eta(p^n,\bm \pi)$ as in (\ref{eq.eta}) and then passing to limit
on $\zz/p^n\zz$. The sequence extends over $R$ (passing to N\'eron
models) and then provides a sequence on component groups
\[0\to \zz \to \phi_{F}\to \qq_p/\zz_p\to 0
\]
over $k$ that coincides with the opposite of the sequence $\zeta_{p^\infty}$ in
(\ref{eq.zeta}). The minus sign depends on Lemma \ref{lem.pullback}.

More generally,  given the $1$-motive
$u_{\pi,K}^2 \colon Y_K\to T_K$ and a character $\chi\in Y^*$,
the sequence $\nu^\vee(\chi)\otimes \zeta$ is obtained as follows:
first consider the push-out $\pre {p^n} T_K\to T_K$ in
$\eta(p^n, u_{\pi,K}^2)$ and then pass to limit on $Y_K/p^nY_K$. At this point
we have a sequence
\begin{equation}\label{eq.limit}
0\to T_K \to F_K\to Y_K\otimes \qq_p/\zz_p\to 0.
\end{equation}
Passing to N\'eron models and taking the induced sequence on component groups
we get a sequence
\begin{equation}\label{eq.comp}
0\to Y^{*\vee} \to \phi_{F}\to Y\otimes \qq_p/\zz_p\to 0.
\end{equation}
This sequence is nothing else than the opposite of push-out with respect to
 $\nu\colon Y\to Y^{*\vee}$ of the sequence $Y\otimes \zeta_{p^\infty}$.
Once fixed a character $\chi\in Y^*$, we can consider the induced
homomorphism $\chi^\vee\colon Y^{*\vee}\to \zz$ (evaluation at
$\chi$). Now the additive extension turns out to be opposite of the push-out
of (\ref{eq.comp}) via the composition of $\chi^\vee$ with the
canonical homomorphism $\zz\to \agr k$.
\item We can also describe
$\nu^\vee(\chi)\otimes \zeta$ in terms of integral of the second
kind generalizing what was done in (\ref{eq.gamma}) and (\ref{eq.h}).

Once fixed
 an isomorphism $Y\cong \oplus_i\zz e_i$, the factor set $\gamma$
in (\ref{eq.gamma})
provides a factor set of $Y\otimes \zeta_{p^\infty}^\add$
\begin{eqnarray*}Y\otimes \gamma \colon ~~  Y\otimes (\qq_p/\zz_p) \times
Y\otimes (\qq_p/\zz_p) &\to &Y\otimes_\zz \zz,\\
\quad \left(\sum_i e_i\otimes a_i, \sum_i e_i\otimes b_i\right)& \mapsto&
\sum_i e_i\otimes [\sigma(a_i)+\sigma(b_i)]
\end{eqnarray*}
and hence a factor set
\begin{eqnarray}\label{fact.add}  Y\otimes (\qq_p/\zz_p) \times
Y\otimes (\qq_p/\zz_p) &\to & \agr R, \nonumber \\
\quad \left(\sum_i e_i\otimes a_i, \sum_i e_i\otimes b_i\right)
&\mapsto&
 \sum_i [\sigma(a_i)+\sigma(b_i)]\mu(e_i,\chi)
\end{eqnarray}
 of $(\nu^\vee(\chi), \id)_* (Y\otimes\zeta^\add_{p^\infty})$.
This factor set  becomes trivial on generic
fibres and a trivialization  is given by
\[h\colon Y_K\otimes \qq_p/\zz_p \to \agr K, \quad \sum_ie_i
\otimes a_i \mapsto \mu(e_i,\chi)\sigma(a_i).\]
If we read this trivialization in terms of formal groups
and then pass to the affine algebras, $h$ corresponds
to a $K$\nobd-homomorphism
\[h^*\colon K[X] \to  K^{\oplus_i\qq_p/\zz_p e_i}, \quad X\mapsto
 \sum_{a\in \oplus_i\qq_p/\zz_p e_i} \sum_i\mu(e_i,\chi)\sigma(a_i)f_a\]
where $f_a\colon \oplus_i\qq_p/\zz_p e_i\to K$ is $1$ in $a$ and
$0$ otherwise.  Also in this case $h^*(X)$
is an integral of the second kind in $K^{\oplus_i\qq_p/\zz_p e_i}$
 that is sent to the class of
 $(\nu^\vee(\chi), \id)_* (Y\otimes\zeta^\add_{p^\infty})$
 via the map \[I_2(Y\otimes \qq_p/\zz_p)\to \MM(Y\otimes \qq_p/\zz_p).\]
\end{enumerate}

\begin{remark} The above constructions  involve the part $u^2_{K, \pi}$
 in Raynaud's decomposition  of the $1-$motive $u_K$. In particular,
the  multiplicative factor set
\begin{eqnarray}\label{fact.mult}
Y_K\otimes (\qq_p/\zz_p) \times Y_K\otimes (\qq_p/\zz_p) &\to&
{\mathbb G}_{m,K}\nonumber \\
\left ( \sum_ie_i \otimes a_i, \sum_ie_i \otimes b_i\right )&\mapsto &
\pi^{-\sum_i[\sigma (a_i)+\sigma (b_i)]\mu(e_i,\chi)},
\end{eqnarray}
whose valuation is the opposite of the additive
factor set in (\ref{fact.add}),
gives  the extension obtained by taking the push-out with respect to
$\chi$ of the  sequence (\ref{eq.limit}).
 If
we think of the valuation as a sort of logarithm killing  elements in $R^*$,
something similar in form, even if quite different in nature,
 happens when working with $1$\nobd-motives of the  type
$u^1_{K, \pi}\colon Y_K\to T_K$, and $u^2_{K, \pi}=0$:
 Also in the present situation  we have a multiplicative factor set,
it may be choosen as follows:
\begin{eqnarray}\label{fact.cl}
\left(\sum_ie_i \otimes a_i, \sum_ie_i \otimes b_i\right )\mapsto \prod_i u(e_i,\chi)^{-[\sigma (a_i)+\sigma (b_i)]},
\end{eqnarray}
where the $u(e_i,\chi)$ are principal units in $R$, so
that the corresponding sequence extends to $R$. Now, we can obtain
from this an additive  factor set by just taking  the $p$-adic logarithm.
As before,
the additive factor set we get,
\[\left ( \sum_ie_i \otimes a_i, \sum_ie_i \otimes b_i\right )
\mapsto \sum_i -[\sigma (a_i)+\sigma (b_i)]\text{log}(u_i)^{}
\] is trivial on generic fibres
and its trivialisation provides an integral of the second kind
$h(\chi)\in K^{\qq_p/\zz_p}$ and  the BT--group of $u^1_{K, \pi}$ is
completely  determined by the
$W(k)$\nobd-module generated  by  the
$h(\chi)$, as $\chi$ varies in the group of the characters  of $T_R$
(cf. \cite{[Fo]}, IV).
Finally, let us observe that if one is just interested in computing
the monodromy, then the use of the valuation is quite appropriate,
but if one  needs to consider
 the integrals of  $u_{K, \pi}=u^1_{K, \pi}+ u^2_{K, \pi}$, {\it i.e.} one
 needs to integrate logarithmic differentials (in the style
of \cite{[Co]}), one  is forced
to extend  the $p-$adic logarithm
 defining  $\log \pi$  in a way allowing one to
distinguish the integrals of $u^1_{K, \pi}$  from those $u^2_{K, \pi}$.
\end{remark}

\section{The finite logarithmic group that extends $\pre n M_K$}

 Throughout this section we work with a strict
$K$-$1$\nobd-motive $u_K\colon \mathbb Z^r\to G_K$ where
$G_K$ is a semiabelian scheme with split
torus part  $\mgr K^d$ and with abelian quotient $A_K$
having good reduction.
Given the $n$-torsion $ \pre {n}M_K $ of such a $1$\nobd-motive,
we know from Theorem \ref{thm.comparison} that it extends to
a finite logarithmic group over $R$,
$ \pre {n}M^{\loga} $, but we still know little concerning
the scheme $\pre {n}M$ underlying
$\pre {n}M^{\loga} $. In this section we are going to describe
the logarithmic group scheme $\pre n M^{\loga} $
 for any $n$ and
precisely we will show that  $\pre n M^{\loga} $ is the valuative space associated to
$(\spec {\pre n A}, {\cal M}_A)$
where $\pre n A$ is  an algebra  over $R$,
 constructed  in a somewhat
 canonical way, and
${\cal M}_A$ is a logarithmic structure on $\spec {\pre n A}$
such that the structure morphism over $R$
 induces a morphism of log schemes
$(\spec {\pre n A},{\cal M}_A)\to \underline T $.

 We have seen that $\pre n M_K$ is an extension of
$Y_K/nY_K$ by $\pre n G$, hence
a $\pre n G_K$\nobd-torsor over $Y_K/nY_K$. If $\pre n G_K={\bm \mu}_n^d$,
 it is easy to describe the algebra of $\pre n M_K$; cfr. \cite{[Mi]},
\cite{[SGA3]}.
For example for the Tate curve $\bm \pi\colon \zz\to \mgr K$
we denoted  $\pre n M_K$ by
\begin{equation}\label{eq.finitek}
\pre n E_K=\spec {  K^{\zz/n\zz}[x]/(x^n-b_{\pi,n}) }
\end{equation}
with $b_{\pi,n}\colon =\sum_{j=0}^{n-1}
\pi^j v_j$,
where $\{v_0, \cdots, v_{n-1}\} $ is the canonical basis of $K^{\zz/n \zz}$.
See also (\ref{eq.finite}).

For the general case, we read from (\ref{dia.HVM}) that
 $\pre n M_K$ is indeed a ${\bm \mu}_n^d$\nobd-torsor
over the finite $K$\nobd-group scheme $\pre n M_K^A$.
Recalling that $\mathrm{Pic}(\pre n M_K^A)=0$,  $\pre n M_K$
can easily be described via  \cite{[Mi]} III \S 4.

\begin{lemma}\label{lem.strM}
Let $ {\cal B}_K$ be the algebra of $\pre n M_K^A$.
Then  $$\pre n M_K=
\spec {{\cal B}_K[T_1,\cdots,T_d,]/(T_1^n-b_1,\dots,T_d^n-b_d)}
$$
for suitable $b_i\in{\cal B}_K^*$, $i=1,\dots,d$.
\end{lemma}

We need now to describe how the $b_i$ above depend on Raynaud's
decomposition of the $1$\nobd-motive $u_K$.

\subsection{The $1$-motive ${\bm \pi}^{-1}$}

Before proceeding with the description of $\pre n M_K$ we need to know
what happens when working with the  $1$-motive
${\bm\pi}^{-1}\colon \zz\to \mgr K, ~1\mapsto \pi^{-1}$.
Its $n$-torsion group scheme is
\begin{equation}\label{eq.finiteminus}
\pre n E^{\frac{1}{\pi}}_K=\spec{\frac{ K^{\zz/n\zz}
[x]}{(x^n-\sum_{i=0}^{n-1}\pi^{-i}v_i)} }
\end{equation}
with $\{v_0,\cdots, v_{n-1}\}$ the canonical base of $ K^{\zz/n\zz}$ over $K$.
It is clear that we can not extend its algebra over $R$
via the equation $x^n-\sum_{i=0}^{n-1}\pi^{-i}v_i$.
However, $\pre n E^{\frac{1}{\pi}}$ is isomorphic to
the group scheme obtained from (\ref{eq.theta})
by push-out with respect to $-1$ or by pull-back with respect
to $-1$. Hence it is also
\begin{equation}\label{eq.finiteminusbis}
\pre n E^{\frac{1}{\pi}}_K\cong \spec{\frac{ K^{\zz/n\zz}
[y]}{(y^n-v_0-\sum_{i=1}^{n-1}\pi^{n-i}v_i)} }
\end{equation}
because  $-1\colon \zz\to \zz$ sends $v_0\mapsto v_0,~ v_i\mapsto
v_{n-i}$ for $i>0$. Moreover this scheme extends over $R$ to
\begin{equation}\label{eq.finiteminusr}
\pre n E^{\frac{1}{\pi}}=
\spec{ \frac{ R^{\zz/n\zz} [y] }{( y^n-v_0-\sum_{i=1}^{n-1}\pi^{n-i}v_i )} }.
\end{equation}
We can endow $\pre n E^{\frac{1}{\pi}}$ with the logarithmic structure
coming from the special fibre and it becomes
 a logarithmic group scheme over $R$. Denote by
$\pre n(\underline{ E}^{\frac{1}{\pi}})$ its valuative logarithmic space;
 this lies in the middle of
$-\theta_n^\pi$
 for $\theta_n^\pi$ as in (\ref{eq.theta}).

Another decomposition of a $1$\nobd-motive that will be useful later
is the following:

\begin{lemma}\label{lem.devplusminus}
Let $u_K$ be a $1$\nobd-motive as in Theorem \ref{thm.devissage}.
Suppose furthermore that $Y_K$ is split and also the torus part
$T_K$ of $G_K$ is split. Once fixed a uniformizing parameter
$\pi\in R$, a basis $(e_j)_j$ of $Y_K\cong \zz^r$ and a basis
$(e_i^*)_i$ of $Y_K^*$  there are decompositions
\[u_K=u_{K,\pi}^1+u_{K,\pi}^2=u_{K,\pi}^1+ u_{K,\pi}^+ + u_{K,\pi}^-
\]
where the first one is the decomposition in Theorem \ref{thm.devissage}. The
second decomposition is uniquely determined by the following conditions:
\begin{itemize} \addtolength{\itemsep}{0.3\baselineskip}
\item $u_{K,\pi}^\pm\colon Y_K\to G_K$ factor through the torus
part.
\item If $\mu^+$ (resp. $\mu^-$) denotes the geometric monodromy
of $ u_{K,\pi}^+$ (resp. of $ u_{K,\pi}^-$), then one has
\mbox{$\mu^+(e_j,e_i^*)\geq 0$} (resp. $\mu(e_j,e_i^*)^-\leq 0$).
\item $u_{K,\pi}^2= u_{K,\pi}^+ + u_{K,\pi}^-$ and $\mu=\mu^+ +
\mu^-$.
\end{itemize}
\end{lemma}
\proof We are reduced to working with $u_{K,\pi}^2$. We define
\[\mu^+ \colon Y_K\otimes Y_K^*\to \zz, \quad
 (e_j,e_i^*)\to
\left\{\begin{array}{lc}\mu(e_j,e_i^*) & \text{ if this is positive} \\
0 & \text{otherwise}
\end{array} \right.\] Similar for $\mu^-$. Moreover $u_{K,\pi}^\pm$ is defined as done for
$u_{K,\pi}^2$ in (\ref{eq.vmotive}) with $\mu^\pm$ in place of $\mu$. \qed

\begin{remark}
Recall the decomposition in Lemma \ref{lem.devplusminus}
and the factorization  $u_{K,\pi}^2=(\nu\otimes \id)\bm \pi_Y$.
Let us construct $\nu^\pm\colon Y_K\to (Y_K^*)^\vee$ from the
geometric monodromy $\mu^\pm$ as in (\ref{eq.nu}).
Then we have factorizations $u_{K,\pi}^+=(\nu^+\otimes \id)\bm \pi_Y$
(resp. $u_{K,\pi}^-=(\nu^-\otimes \id)\bm \pi^{-1}_Y$)
with
$$\bm \pi^{-1}_Y\colon Y_K\to Y_K\otimes_\zz\mgr K, ~~
y\mapsto y\otimes \pi^{-1} .$$
\end{remark}

\subsection{A general result on short exact sequences}\label{sect.vert}
Looking at the diagram (\ref{dia.HVM}) where $\pre n T_K={\bm
\mu}_n^d$ we realize that we should move our attention from the
horizontal sequence in the middle $\eta(n,u_K)$ to the vertical
sequence in the middle, because, as we have already observed ${\bm
\mu}_n$-torsors can easily be described. Moreover we need to know
how such a "vertical sequence" depends on the analogous sequences
for $u_{\pi,K}^1$ and $u_{\pi,K}^2$. For this we need a general
result.

Let $\bar\psi\colon 0 \to I \overset{w}{\to} L\overset{k}{\to} P
\to 0$ be an exact sequence of group schemes over a base scheme
$S$.  Let $N$ be another group scheme over
$S$ and consider two extensions $\eta^i\colon 0\to L\to M^i\to
N\to 0$, $i=1,2$. Let $\eta\colon 0\to L\to M\to N\to 0$ be a
sequence isomorphic to the Baer sum $\eta^1+\eta^2$. Consider
then the following diagrams
\[\xymatrix{  & {I} \ar@{^{(}->}[d]^w \ar@2{-}[r]& I \ar@{^{(}->}[d]^{\tau^i} & &  \\
\eta^i\colon 0 \ar[r] &L \ar@{->>}[d]^k \ar[r] &
              M^i  \ar@{->>}[d]^{g^i} \ar[r]^{h^i} & N \ar@2{-}[d] \ar[r] & 0\\
~~~~~0 \ar[r] & P  \ar[r] &
              Q^i \ar[r]^{f^i} & N\ar[r] & 0
}\quad\quad
\xymatrix{ & {I} \ar@{^{(}->}[d]^w \ar@2{-}[r]& I \ar@{^{(}->}[d]^{\tau} & &  \\
\eta\colon 0 \ar[r] &L \ar@{->>}[d]^k \ar[r] &
              M  \ar@{->>}[d]^{g} \ar[r]^{h} & N \ar@2{-}[d] \ar[r] & 0\\
~~~~~0 \ar[r] & P  \ar[r] &
              Q \ar[r]^{f} & N\ar[r] & 0
}
\]
where the vertical sequence on the left is $\bar\psi$ and the upper horizontal
sequence is $\eta^i$ for $i=1,2$ (resp. $\eta$).
Call $\psi^i$ for $i=1,2$ (resp. $\psi$) the vertical sequence in the middle.
Suppose now that there is a sequence
\[
\tilde\eta^2 :\quad 0\to I\to \tilde M^2\to N\to 0
\]
such that $\eta^2=w_*\tilde \eta^2$.
 Summarizing we have
\[\framebox{ $\eta \cong \eta^1+w_*\tilde \eta^2$. } \]
We are going to see that a similar relation holds also for the vertical
sequences $\psi,\psi^1$.

Now $k_*w_*\tilde\eta^2$ is isomorphic to the trivial extension
and we choose a section $\sigma$ of $f^2$. Moreover
  $k_*\eta\cong k_*\eta^1$ and there is then an isomorphism $
{\iota}^\sigma\colon Q \to Q^1$ depending on $\sigma$.
It is also not difficult to check that
\begin{equation}\label{eq.isopsigeneral}
\psi\cong (\iota^\sigma)^*\psi^1+ (\sigma f)^*\psi^2
\end{equation}
Consider also the following   push-out diagram
\[\xymatrix{ \tilde \eta^2\colon & 0 \ar[r] & I
 \ar[d]^w \ar@{-->}[dr]^{\tau^2} \ar[r]& \tilde M^2 \ar[d]^\delta
\ar[r]^{\tilde h^2} & N \ar@2{-}[d]
 \ar[r] & 0 \\
w_*\tilde \eta^2=\eta^2\colon& 0 \ar[r] &L \ar[d]^k \ar[r] &
              M^2 \ar[d]^{g^2} \ar[r]^{h^2} & N \ar@2{-}[d] \ar[r] & 0\\
k_*\eta^2\colon & 0 \ar[r] & P  \ar[r] &
              Q^2 \ar[r]^{f^2} & N\ar[r] & 0
}
\]
where lower sequence is isomorphic to the trivial one.
We denoted by  $\psi^2$  the exact sequence
involving $g^2$ and $\tau^2$.
There exists then a canonical homomorphism
$\sigma^c\colon N\to Q^2$ such that
\[\sigma^c \tilde h^2=g^2 \delta ; \]
 it is not difficult to check that $\sigma^c$ is a section of $f^2$.
Moreover by construction it satisfies
\[(\sigma^{c})^*\psi^2\cong \tilde \eta^2\]
Taking now $\sigma=\sigma^c$ in (\ref{eq.isopsigeneral}) and setting
 $\iota\colon=\iota^{\sigma^c} $ we get that
\begin{equation}\label{eq.psi}
\framebox{ $\psi\cong \iota^*\psi^1+ f^*\tilde\eta^2 $ .}
\end{equation}

\subsection{${}_n M_K$ as torsor under ${\bm \mu}_n^d$}

 Let $u_K$ be a $1$\nobd-motive as in Theorem \ref{thm.devissage}.
We have then a decomposition $u_K=u_{K,\pi}^1+u_{K,\pi}^2$
and an isomorphism of sequences
$\eta(n,u_K)\cong \eta(n,u_{K,\pi}^1)
+\eta(n,u_{K,\pi}^2)$ for $\eta(n,-)$ the sequence introduced in
(\ref{eq.eta}).

The $n$-torsion $ \pre {p^n}M_K$ of $u_K$ lies
in the middle of $\eta(n,u_K)$
and we have already seen what it looks like
in Lemma~\ref{lem.strM}.
Moreover $\eta(n,u_{K,\pi}^2)=w_*\eta(n,\tilde u_K^2)$,
where the $1$\nobd-motive
$\tilde u_K^2\colon Y_K\to T_K$ is obtained from $u_{K,\pi}^2$
by forgetting the inclusion $T_K\to G_K$ and $w$ is this inclusion
restricted to kernels of $n$\nobd-multiplication.
Hence we are in the situation of the subsection \ref{sect.vert}.

Considering diagram (\ref{dia.HVM}) for $\eta(n,u_{K,\pi}^1)$
in place of $\eta(n,u_{K,\pi} )$ we get
\begin{equation}\label{dia.HVM1}
\xymatrix{ &   & \pre n T_K \ar[d]^w \ar@2{-}[r]& \pre n T_K \ar[d]^{\tau^1} \\
 \eta(n,u_{K,\pi}^1) \colon &0 \ar[r] &\pre n G_K \ar[d]^{\rho_n} \ar[r] &
              \pre n M_K^1 \ar[d]^{g^1} \ar[r]^{h^1} & Y_K/nY_K \ar@2{-}[d] \ar[r] & 0\\
\eta(n,\rho u_{K,\pi}^1)\colon &0 \ar[r] &\pre n A_K  \ar[r] &
              Q^1_K \ar[r]^{f^1} & Y_K/nY_K \ar[r] & 0   .
}
\end{equation}
Let $\psi(n,u_{K,\pi}^1)$ be the vertical sequence in the middle
of this diagram and $\psi(n,u_{K,\pi} )$ the corrisponding vertical
sequence in the middle of  (\ref{dia.HVM}). By (\ref{eq.psi}) we have
\begin{equation}\label{iso.psimot}
\psi(n,u_K)\cong \iota^*\psi(n,u_{K,\pi}^1)
+f^*\eta(n,\tilde u_K^2)
\end{equation}
for $\iota$ the isomorphism $\pre n M_K^A\to Q_K^1$ constructed
as in \S\ref{sect.vert}
By abuse of notations, we denote by  $\pre n M_K^1$ also the group scheme
 in the middle of $ \iota^*\psi(n,u_{K,\pi}^1)$; indeed it
is isomorphic to the one in the middle of $\eta(n,u_{K,\pi}^1)$. The scheme
 $\pre n M_K^1$ is a $\bm \mu_n^d$\nobd-torsor
over $\pre n M_K^A$ and hence
\[\pre n M_K^1 \cong
\spec {{\cal B}_K[T_1,\cdots,T_d]/(T_1^n-b_1^{(1)}, \dots, T_d^n-b_d^{(1)}) }. \]
On the other hand, call $\pre n \tilde M_K^2$ the group in the
middle of $\eta(n,\tilde u_K^2)$. It is a $\bm \mu_n^d$\nobd-torsor
over $\zz^r/n\zz^r$ and hence it  has the form
\[\pre n \tilde M_K^2 \cong
\spec {K^{\zz^r/n\zz^r}[T_1,\cdots,T_d]/(T_1^n-b_1^{(2)},\dots,T_d^n-b_d^{(2)}) }.
\]
Hence the group
in the middle of $f^* \eta(n,\tilde u_{K,\pi}^2)$  will be

\[ \spec { {\cal B}_K  [T_1,\cdots,T_d]/(T_1^n-b_1^{(2)},\dots,T_d^n-b_d^{(2)} ) } .
\]
Recall now that $ \tilde u_{K,\pi}^2=u_{K,\pi}^+ + u_{K,\pi}^-$
by Lemma \ref{lem.devplusminus} and hence
$\eta(n, u_{K,\pi}^2)\cong\eta(n, u_{K,\pi}^+)+\eta(n, u_{K,\pi}^-)$.
Let \[ \spec { K^{\zz^r/n\zz^r}  [T_1,\cdots,T_d]/(T_1^n-b_1^+,\dots,T_d^n-b_d^+) }
\]
be the group scheme in the middle of $\eta(n,u_{K,\pi}^+)$
and analogously for $u_{K,\pi}^-$ with elements
\mbox{$b_i^-\in K^{\zz^r/n\zz^r}$} in place of $b_i^+$.
Using the sequence
\begin{equation}\label{seq.milnetorsor}
 0\to \bm\mu_n^d(X) \to \Gamma(X,{\cal O}_X)^{d*} \stackrel{n}{\to} \Gamma(X,{\cal
O}_X)^{d*} \to {\mathrm H}(X,\bm \mu_n^d)\to 0
\end{equation}
we may assume that $b_i^{(2)}=b_i^+b_i^-$
and
\[ b_i =b_i^{(1)}b_i^+b_i^- ~~\text{for all}~~ 1\leq i\leq d . \]

Recall now that the vertical sequence on the left in (\ref{dia.HVM})
or  (\ref{dia.HVM1})
extends over $R$ because $G_K$ has good reduction and that
also $\pre n M^A_K$ extends to a finite group scheme over $R$, say $\pre n
M^A=\spec {\cal B}$, because the $K$\nobd-$1$\nobd-motive
$$\rho u_K\colon \zz^r\to G_K\to A_K$$
 has good
reduction. Finally also $\eta(n,u_{K,\pi}^1)$ extends over $R$ and
hence the same is true for the vertical sequence
$\psi(n,u_{K,\pi}^1) $ and so we may assume $b_i^{(1)}\in {\cal
B}^*$. This implies that if we want  to extend $\pre n M_K$ over
$R$ we have to understand better $\pre n\tilde  M_K^2$ and hence
$b_i^{(2)}=b_i^+b_i^-$. Recall the description of $u_{K,\pi}^2$ in
(\ref{eq.vmotive}) where now $Y_K=\zz^r$ and $T_K=\mgr K^d$. It is
an easy exercise to check that
\[
b_i^+= (b_{\pi,n})^{ \sum_{j=1}^r{\mu^+(e_j,e_i^*)} }
~~\text{with}~~ b_{\pi,n}\colon =\sum_{j=0}^{n-1} \pi^j v_j\quad
\in R^{\zz^r/n \zz^r}, \] with $v_0,\dots,v_{n-1}$  the standard
basis of $K^{\zz/n \zz}$ (resp. of $R^{\zz/n \zz}$),
$e_1,\dots,e_r$ the usual basis of $\zz^r$ and $e_i^*$ the character
in $Y^*_K$ such that $e_i^*(T_h)=\delta_{ih}$, $\mu^+$ the
geometric monodromy of $u_{K,\pi}^+$. Moreover
$\mu^+(e_j,e_i^*)\geq 0$ by definition of $u_{K,\pi}^+$; hence
$b_i^+\in {\cal B}$ for all $i$. In a similar way
\[\displaystyle{b_i^-  = (b_{\pi^{-1},n})^{  \sum_{j=1}^r {-\mu^-(e_j,e_i^*)} }
~~\text{with}~~ b_{\pi^{-1},n}\colon =
v_0+\sum_{i=1}^{n-1}\pi^{n-i}v_i \in  R^{\zz/n\zz}}  .
\]
Hence also $b_i^-\in {\cal B}$ for all $i$ because $\mu^-(e_j,e_i^*)
\leq 0$ by definition of $u_{K,\pi}^-$. Summarizing,
$b_i=b_i^{(1)}b_i^{(2)}=b_i^{(1)}b_i^+b_i^-\in {\cal B}$ for all
$i$. Hence $\pre n M_K$ extends to a finite scheme over $\spec
{\cal B}$
\begin{equation}\label{eq.mscheme}
\pre n M =  \spec {{\cal B}[T_1,\cdots,T_d]/(T_1^n-b_1,\dots, T_d^n-b_d)}.
 \end{equation}
The following lemma says that $\pre n M$ is  indeed a ``nice''
model.

\begin{lemma}
Let notations be as above and endow $\pre n M^A=\spec {\cal B}$
 with the inverse image log structure of the base $\underline T$.
Let  $(\pre n M,{\cal M}_{\pre n M} )$ be the
scheme (\ref{eq.mscheme}) with the logarithmic structure induced by the
special fibre. Then the canonical morphism
$\pre n M_K\to \pre n M^A_K $
extends to a unique morphism of logarithmic groups over $R$.
Moreover the following universal property holds:
For any   fine saturated logarithmic scheme
$(S,{\cal M}_S)$ over $(\pre n M^A,{\cal M}_{\pre n M^A})$
there is a bijection between the scheme theoretic morphisms
$S_K\to \pre n M_K$ over $\pre n M_K^A$ and the
logarithmic morphisms $(S,{\cal M}_S)\to (\pre n M,{\cal M}_{\pre n M} )$
over $(\pre n M^A,{\cal M}_{\pre n M^A})$.
\end{lemma}
\proof
We may assume $S=\spec C$ affine.
Any $\pre n M^A_K$\nobd-morphism $S_K\to \pre n M_K$
is described as a homomorphism of ${\cal B}_K$\nobd-algebras
$\varphi_K\colon {\cal B}_K[T_1,\cdots,T_n]/(T_i^n-b_i)\to C_K$.
The only problem for the extension
is to prove that the images of all $T_i$ lie in $C$ and more precisely in
$\Gamma(S,{\cal M}_S)$. However  $b_i=b_i^{(1)}b_i^{(2)}$ with
 $b_i^{(1)}\in  {\cal B}^*$ and $b_i^{(2)}\in \Gamma(Y,{\cal M}_Y)$
for any logarithmic space over the group $\zz/n\zz$ endowed with
the inverse image log structure because of the description of
$b_i^{(2)}$ in terms of $b_{\pi,n}$  and
$b_{\pi^{-1},n}$. Moreover $\varphi_K(T^n_i)=\varphi_K(b_i)\in
\Gamma(S,{\cal M}_S )$. It is now sufficient to recall that $S$ is
saturated to conclude that also $\varphi_K(T_i)=b_i\in
\Gamma(S,{\cal M}_S )$. \qed

\begin{proposition}
Let notations be as above. The finite group scheme
$\pre n M_K$ over $\pre n M_K^A$
 extends (up to isomorphisms) to a finite logarithmic group $\pre n M^\loga$
that is the valuative logarithmic space associated to a logarithmic
scheme whose  underlying scheme is
\begin{equation}\label{eq.mlog}
\pre n M=\spec {{\cal B}[T_1,\cdots,T_d]/(T_1^n-b_1,\dots, T_d^n-b_d) }.
 \end{equation}
Moreover the diagram (\ref{dia.HVM}) extends to a diagram
of finite logarithmic groups with $\pre n M^\loga$ in the middle.
\end{proposition}
\proof
It remains only to prove that $\pre n M$ with the logarithmic structure
induced by the special fibre induces a group functor on the category
of fine saturated logarithmic schemes over $R$, i.e. it is in $T^\loga_\fl$.
However this is immediate consequence of the previous lemma.

Also the assertion on the diagram follows applying  the previous lemma.
\qed
\medskip

Observe that we had already proved in Theorem
\ref{prop.comparison} that $\pre n M_K$ extends to a logarithmic
group over $\underline T$, however, in the hypothesis of this
section, it is possible to describe it in terms of algebras and not
only as a ``sum'' of two extensions.

\end{document}